\documentclass[11pt]{article}
\usepackage{amssymb}
\newtheorem{lemma}{Lemma}[section]
\newtheorem{definition}[lemma]{Definition}

\newtheorem{theorem}[lemma]{Theorem}
\newtheorem{proposition}[lemma]{Proposition}
\newtheorem{corollary}[lemma]{Corollary}
\newtheorem{example}[lemma]{Example}

\newtheorem{note}[lemma]{Remark}

\def\endproof{\hfill$\Box$}

\def\dcl{\mathop{\rm dcl}\nolimits}

\def\endproof{\hfill$\Box$}

\title{On relative separability \\ in hypergraphs of models of theories}

\author{B.Sh.~Kulpeshov, S.V.~Sudoplatov}
\date{}

\begin{document}

\maketitle

\begin{abstract}
In the paper, notions of relative separability for hypergraphs of
models of a theory are defined. Properties of these notions and
applications to ordered theories are studied: characterizations of
relative separability both in a general case and for almost
$\omega$-ca\-te\-go\-ri\-cal quite o-minimal theories are
established.
\end{abstract}


\bigskip
Hypergraphs of models of a theory are related to derivative
objects allowing to obtain an essential structural information
both on theories themselves and related semantical objects
including graph ones  \cite{CCMCT14, Su013, Su08, Baik, SudKar16,
KulSud17, Sud17, KulSud171}.

In the present paper, notions of relative separability for
hypergraphs of models of a theory are defined. Properties of these
notions and applications to ordered theories are studied:
characterizations of relative separability both in a general case
and for almost $\omega$-categorical quite o-minimal theories are
established.

\section{Preliminaries}
\noindent

Recall that a {\it hypergraph\/} is a pair of sets $(X,Y)$, where
$Y$ is some subset of the Boolean $\mathcal{P}(X)$ of the set $X$.

Let $\mathcal{M}$ be some model of a complete theory $T$.
Following \cite{SudKar16}, we denote by $H(\mathcal{M})$ a family
of all subsets $N$ of the universe $M$ of $\mathcal{M}$ that are
universes of elementary submodels $\mathcal{N}$ of the model
$\mathcal{M}$: $H(\mathcal{M})=\{N\mid
\mathcal{N}\preccurlyeq\mathcal{M}\}$. The pair
$(M,H(\mathcal{M}))$ is called the {\it hypergraph of elementary
submodels\/} of the model $\mathcal{M}$ and denoted by
$\mathcal{H}(\mathcal{M})$\index{$\mathcal{H}(\mathcal{M})$}.

For a cardinality $\lambda$ by
$H_\lambda(\mathcal{M})$\index{$H_\lambda(\mathcal{M})$} and
$\mathcal{H}_\lambda(\mathcal{M})$\index{$\mathcal{H}_\lambda(\mathcal{M})$}
are denoted restrictions for $H(\mathcal{M})$ and
$\mathcal{H}(\mathcal{M})$ respectively on the class of elementary
submodels $\mathcal{N}$ of models $\mathcal{M}$ such that
$|N|<$~$\lambda$.

By
$\mathcal{H}_p(\mathcal{M})$\index{$\mathcal{H}_p(\mathcal{M})$}
we denote the restriction of the hypergraph
$\mathcal{H}_{\omega_1}(\mathcal{M})$ on the class of elementary
submodels ${\cal N}$ of the model ${\cal M}$ that are prime over
finite sets. Similarly by
$H_p(\mathcal{M})$,\index{$H_p(\mathcal{M})$} is denoted the
corresponding restriction for $H_{\omega_1}(\mathcal{M})$.

\begin{definition} {\rm \cite{SudKar16, Engel}.
Let $(X,Y)$ be a hypergraph, $x_1,x_2$ be distinct elements of
$X$. We say that the element $x_1$ is {\em separated} or {\em
separable} from the element $x_2$, or {\em $T_0$-separable} if
there is $y\in Y$ such that $x_1\in y$ and $x_2\notin y$. The
elements $x_1$ and $x_2$ are called {\em separable}, {\em
$T_2$-separable}, or {\em Hausdorff separable} if there are
disjoint $y_1,y_2\in Y$ such that $x_1\in y_1$ and $x_2\in y_2$.}
\end{definition}

\begin{theorem}\label{th21} {\rm \cite{SudKar16}.} Let $\mathcal{M}$ be an $\omega$-saturated model of
a countable complete theory $T$, $a$ and $b$ be elements of $\mathcal{M}$. The following are equivalent:

$(1)$ the element $a$ is separable from the element $b$ in
$\mathcal{H}(\mathcal{M})$;

$(2)$ the element $a$ is separable from the element $b$ in
$\mathcal{H}_{\omega_1}(\mathcal{M})$;

$(3)$ $b\notin{\rm acl}(a)$.
\end{theorem}

\begin{theorem}\label{th26} {\rm \cite{SudKar16}.} Let $\mathcal{M}$ be an $\omega$-saturated model of
a countable complete theory $T$, $a$ and $b$ be elements of $\mathcal{M}$. The following are equivalent:

$(1)$ the elements $a$ and $b$ are separable in
$\mathcal{H}(\mathcal{M})$;

$(2)$ the elements $a$ and $b$ are separable in
$\mathcal{H}_{\omega_1}(\mathcal{M})$;

$(3)$ ${\rm acl}(a)\cap{\rm acl}(b)=\varnothing$.
\end{theorem}

\begin{corollary}\label{co27} {\rm \cite{SudKar16}.}
Let $\mathcal{M}$ be an $\omega$-saturated model of a countable complete theory $T$, $a$ and $b$ be elements
of $\mathcal{M}$, and there exists the prime model over $a$. The following are equivalent:

$(1)$ the element $a$ is separable from the element $b$ in
$\mathcal{H}(\mathcal{M})$;

$(2)$ the element $a$ is separable from the element $b$ in
$\mathcal{H}_{\omega_1}(\mathcal{M})$;

$(3)$ the element $a$ is separable from the element $b$ in
$\mathcal{H}_p(\mathcal{M})$;

$(4)$ $b\notin{\rm acl}(a)$.
\end{corollary}

\begin{corollary}\label{co28} {\rm \cite{SudKar16}.}
Let $\mathcal{M}$ be an $\omega$-saturated model of a countable complete theory $T$, $a$ and $b$ be elements
of $\mathcal{M}$, and there exist the prime models over $a$ and $b$ respectively. The following are
equivalent:

$(1)$ the elements $a$ and $b$ are separable in
$\mathcal{H}(\mathcal{M})$;

$(2)$ the elements $a$ and $b$ are separable in
$\mathcal{H}_{\omega_1}(\mathcal{M})$;

$(3)$ the elements $a$ and $b$ are separable in
$\mathcal{H}_p(\mathcal{M})$;

$(4)$ ${\rm acl}(a)\cap{\rm acl}(b)=\varnothing$.
\end{corollary}

\begin{definition} {\rm \cite{SudKar16}.
Let $(X,Y)$ be a hypergraph, $X_1,X_2$ be disjoint nonempty
subsets of the set $X$. We say that the set $X_1$ is {\em
separated} or {\em separable} from the set $X_2$, or {\em
$T_0$-separable} if there is $y\in Y$ such that $X_1\subseteq y$
and $X_2\cap y=\varnothing$. The sets $X_1$ and $X_2$ are called
{\em separable}, {\em $T_2$-separable}, or {\em Hausdorff
separable} if there are disjunct $y_1,y_2\in Y$ such that
$X_1\subseteq y_1$ and $X_2\subseteq y_2$.}
\end{definition}

By using proofs of theorems \ref{th21} and \ref{th26}, the
following generalizations of these theorems are established.

\begin{theorem}\label{th31} {\rm \cite{SudKar16}} Let $\mathcal{M}$
be a $\lambda$-saturated model of a complete theory $T$,
$\lambda\geq{\rm max}\{|\Sigma(T)|$, $\omega\}$, $A$ and $B$ be
nonempty sets in $\mathcal{M}$ having the cardinalities
$<\lambda$. The following are equivalent:

$(1)$ the set $A$ is separable from the set $B$ in
$\mathcal{H}(\mathcal{M})$;

$(2)$ the set $A$ is separable from the set $B$ in
$\mathcal{H}_\lambda(\mathcal{M})$;

$(3)$ ${\rm acl}(A)\cap B=\varnothing$.
\end{theorem}

\begin{theorem}\label{th32} {\rm \cite{SudKar16}} Let $\mathcal{M}$ be a $\lambda$-saturated model of a complete
theory $T$,  $\lambda\geq{\rm max}\{|\Sigma(T)|$, $\omega\}$, $A$
è $B$ be nonempty sets in $\mathcal{M}$ having the cardinalities
$<\lambda$. The following are equivalent:

$(1)$ the sets $A$ and $B$ are separable in
$\mathcal{H}(\mathcal{M})$;

$(2)$ the sets $A$ and $B$ are separable in
$\mathcal{H}_\lambda(\mathcal{M})$;

$(3)$ ${\rm acl}(A)\cap{\rm acl}(B)=\varnothing$.
\end{theorem}

We obtain by analogy with corollaries \ref{co27} and \ref{co28}

\begin{corollary}\label{co33} {\rm \cite{SudKar16}.}
Let $\mathcal{M}$ be an $\omega$-saturated model of a small theory
$T$, $A$ and $B$ be finite nonempty sets in $\mathcal{M}$. The
following are equivalent:

$(1)$ the set $A$ is separable from the set $B$ in
$\mathcal{H}(\mathcal{M})$;

$(2)$ the set $A$ is separable from the set $B$ in
$\mathcal{H}_{\omega_1}(\mathcal{M})$;

$(3)$ the set $A$ is separable from the set $B$ in
$\mathcal{H}_p(\mathcal{M})$;

$(4)$ ${\rm acl}(A)\cap B=\varnothing$.
\end{corollary}

\begin{corollary}\label{co34} {\rm \cite{SudKar16}.}
Let $\mathcal{M}$ be an $\omega$-saturated model of a small theory
$T$, $A$ and $B$ be finite nonempty sets in $\mathcal{M}$. The
following are equivalent:

$(1)$ the sets $A$ and $B$ are separable in
$\mathcal{H}(\mathcal{M})$;

$(2)$ the sets $A$ and $B$ are separable in
$\mathcal{H}_{\omega_1}(\mathcal{M})$;

$(3)$ the sets $A$ and $B$ are separable in
$\mathcal{H}_p(\mathcal{M})$;

$(4)$ ${\rm acl}(A)\cap{\rm acl}(B)=\varnothing$.
\end{corollary}

The following proposition extends Theorem \ref{th32} with an
additional criterion.

\begin{proposition}\label{genpr_asets}
Let $T$ be a theory, $\mathcal{M}\models T$, $\emptyset\ne
A\subseteq M$, $\emptyset\ne B\subseteq M$, $\mathcal{M}$ be
$|A\cup B|^+$-saturated. Then $A$ and $B$ are separable from each
other in $\mathcal{H}(\mathcal{M})$ if and only if the following
conditions hold:

$(1)$ ${\rm acl}(A)\cap{\rm acl}(B)=\emptyset$;

$(2)$ For any isolated type $p\in S_1(\emptyset)$,
$p(\mathcal{M})\setminus {\rm acl}(A)\ne \emptyset$ and
$p(\mathcal{M})\setminus {\rm acl}(B)\ne \emptyset$.
\end{proposition}

Proof of Proposition \ref{genpr_asets}. If $A$ and $B$ are
separable from each other in $\mathcal{H}(\mathcal{M})$ then by
Theorem \ref{th32} we have ${\rm acl}(A)\cap{\rm
acl}(B)=\emptyset$. If there is an isolated type $p\in
S_1(\emptyset)$ such that $p(M)\subseteq {\rm acl}(A)$ then there
is $\mathcal{M}_2\prec \mathcal{M}$ with $B\subseteq M_2$ and
$p(\mathcal{M})\cap M_2=\emptyset$, i.e. $p$ is not realized in
$\mathcal{M}_2$. Similarly, $p(\mathcal{M})\not\subseteq {\rm
acl}(A)$.

If the conditions (1), (2) hold then $A$ and $B$ are separable
from each other in $\mathcal{H}(\mathcal{M})$ by Theorem
\ref{th32}.
\endproof

Recall that a subset $A$ of a linearly ordered structure $M$ is called {\it convex}
if for any $a, b\in A$ and $c\in M$ whenever $a<c<b$ we have $c\in A$. A {\it weakly
o-minimal structure} is a linearly ordered structure $M=\langle M,=,<,\ldots \rangle$
such that any definable (with parameters) subset of the structure $M$ is a union of
finitely many convex sets in $M$.

In the following definitions $M$ is a weakly o-minimal structure, $A, B\subseteq M$,
$M$ be $|A|^+$-saturated, $p,q\in S_1(A)$ be non-algebraic types.
\begin{definition}\rm \cite{bbs1}
We say that $p$ is not {\it weakly orthogonal} to $q$  ($p\not\perp^w  q$) if there exist
an $A$-definable formula $H(x,y)$, $\alpha \in p(M)$ and $\beta_1, \beta_2 \in q(M)$
such that $\beta_1 \in H(M,\alpha)$ and $\beta_2 \not\in H(M,\alpha)$.
\end{definition}

\begin{definition} \rm \cite{k2003} We say that $p$ is not {\it quite orthogonal} to
 $q$ ($p\not\perp^q q$) if there exists an $A$-definable bijection $f: p(M)\to q(M)$.
 We say that a weakly o-minimal theory is {\it quite o-minimal} if the notions of
 weak and quite orthogonality of 1-types coincide.
\end{definition}

In the work \cite{KS} the countable spectrum for quite o-minimal theories with non-maximal
number of countable models has been described:

\begin{theorem}\label{KS_apal} Let $T$ be a quite o-minimal theory with non-maximal number
of countable models. Then $T$ has exactly $3^k\cdot 6^s$ countable models, where $k$ and
$s$ are natural numbers. Moreover, for any $k,s\in\omega$ there exists a quite o-minimal
theory $T$ having exactly $3^k\cdot 6^s$ countable models.
\end{theorem}

Realizations of these theories with a finite number of countable models are natural
ge\-ne\-ra\-li\-za\-tions of Ehrenfeucht examples obtained by expansions of dense linear orderings by
a countable set of constants, and they are called theories of  {\em Ehrenfeucht type}. Moreover,
these realizations are representative examples for hypergraphs of prime models
\cite{CCMCT14, Su08, SudKar16}.

\section{Relative separability in hypergraphs of models of theories}

Observe that since by Theorem \ref{th32} and Corollary \ref{co34}
separability of sets $A$ and $B$ in hypergraphs $H(\mathcal{M})$
is possible only when ${\rm acl}(A)\cap{\rm acl}(B)=\varnothing$,
such a separability doesn't hold when ${\rm
acl}(\emptyset)\ne\emptyset$. Thus, it is natural to consider the
following notions of {\em relative} separability.

\begin{definition} {\rm
Let $(X,Y)$ be a hypergraph, $x_1,x_2$ be distinct elements of
$X$, $Z\subset X$, $x_2\notin Z$. We say that the element $x_1$ is
{\em $Z$-separated} or {\em $Z$-separable} from the element $x_2$,
or {\em $(T_0,Z)$-separable} if there is $y\in Y$ such that
$x_1\in y\cup Z$ and $x_2\notin y$. In this case the set $y$ is
called {\em $Z$-separating} $x_1$ from $x_2$. At the additional
condition $x_1\notin Z$ the elements $x_1$ and $x_2$ are called
{\em $Z$-separable}, {\em $(T_2,Z)$-separable}, or {\em Hausdorff
$Z$-separable} if there are $y_1,y_2\in Y$ such that $(y_1\cap
y_2)\setminus Z=\varnothing$, $x_1\in y_1$ and $x_2\in y_2$.

Let $X_1,X_2$ be nonempty subsets of the set $X$, $(X_1\cap
X_2)\setminus Z=\emptyset$, $X_2\not\subseteq Z$. We say that the
set $X_1$ is {\em $Z$-separated} or {\em $Z$-separable} from the
set $X_2$, or {\em $(T_0,Z)$-separable} if there is $y\in Y$ such
that $X_1\subseteq y\cup Z$ and $(X_2\cap y)\setminus
Z=\varnothing$. At the additional condition $X_1\not\subseteq Z$
the sets $X_1$ and $X_2$ are called {\em $Z$-separable}, {\em
$(T_2,Z)$-separable}, or {\em Hausdorff $Z$-separable} if there
are $y_1,y_2\in Y$ such that $(y_1\cap y_2)\setminus
Z=\varnothing$, $X_1\subseteq y_1\cup Z$ and $X_2\subseteq y_2\cup
Z$.}
\end{definition}

\begin{note}\label{rem220}
{\rm 1. The notions of separability given in Section 1 correspond $Z$-separability for
$Z=\varnothing$, $X_1\ne\varnothing$, $X_2\ne\varnothing$.

2. If $X_2\subseteq Z$ then the set $X_2$ also can be assumed
$Z$-separable from $X_1$, although there is no reason to say on real separability
of elements of the set $X_2$ from $X_1$.}
\end{note}

For a tuple $\bar{a}$ and a set $Z$ we denote by $\bar{a}Z$ the union of the set $Z$
with the set of all elements containing in $\bar{a}$.

The following theorem modifies Theorem \ref{th21}, and it is a generalization of the theorem for ${\rm
acl}(\varnothing)=\varnothing$.

\begin{theorem}\label{th221} Let $\mathcal{M}$ be an $\omega$-saturated
model of a countable complete theory $T$, $Z$ be the algebraic
closure of some finite set in $\mathcal{M}$, $a$ and $b$ be
elements of $\mathcal{M}$, $b\notin Z$. The following are
equivalent:

$(1)$ the element $a$ is $Z$-separable from the element $b$ in $H(\mathcal{M})$
by some set $y$ from $H(\mathcal{M})$ containing $Z$;

$(2)$ the element $a$ is $Z$-separable from the element $b$ in $H_{\omega_1}(\mathcal{M})$
by some set $y$ from $H_{\omega_1}(\mathcal{M})$ containing $Z$;

$(3)$ $b\notin{\rm acl}(aZ)$.
\end{theorem}

Proof. The implications $(2)\Rightarrow(1)$ and $(1)\Rightarrow(3)$ are obvious (clearly, if
$b\in{\rm acl}(Z\cup\{a\})$ then $b$ belongs to any model $\mathcal{N}\preccurlyeq\mathcal{M}$
containing $Z\cup\{a\}$).

To prove the implication $(3)\Rightarrow(2)$ we need the following lemma.

\begin{lemma}\label{le224}
Let $\bar{a}$ be a tuple, $B$ be a finite set for which $({\rm
acl}(\bar{a}Z)\cap B)\setminus Z=\varnothing$, and
$\varphi(x,\bar{a})$ be some consistent formula. Then there is an
element $c\in\varphi(\mathcal{M},\bar{a})$ such that $({\rm
acl}(\bar{a}cZ)\cap B)\setminus Z=\varnothing$.
\end{lemma}

Proof of Lemma \ref{le224}. If $\varphi(\mathcal{M},\bar{a})\cap Z\ne\varnothing$ then there is nothing
to prove since as $c$ we can take an arbitrary element of $\varphi(\mathcal{M},\bar{a})\cap Z$.

Suppose that $\varphi(\mathcal{M},\bar{a})\cap Z=\varnothing$.
By compactness and using consistent formulas $\varphi'(x,\bar{a})$ with the condition
$\varphi'(x,\bar{a})\vdash\varphi(x,\bar{a})$ instead of $\varphi(x,\bar{a})$, it is sufficiently
to prove that for any $d\in B\setminus Z$ and the finite set of formulas
$\psi_1(x,\bar{a},y),\ldots,\psi_n(x,\bar{a},y)$ with the condition
$$\psi_i(x,\bar{a},y)\vdash\varphi'(x,\bar{a})\wedge\forall x\left(\varphi'(x,\bar{a})\to\exists^{=k_i}y\psi_i(x,\bar{a},y)\right)$$
for some natural $k_i$, $i=1,\ldots,n$, there is an element $c\in\varphi'(\mathcal{M},\bar{a})$
such that $$\models\bigwedge\limits_{i=1}^n\neg\psi_i(c,\bar{a},d).$$

Assume the contrary: for any $c\in\varphi'(\mathcal{M},\bar{a})$ there is $i$ such that
$\models\psi_i(c,\bar{a},d)$. Then the formula
$\chi(x,\bar{a},y)\rightleftharpoons\bigvee\limits_{i=1}^n\psi_i(x,\bar{a},y)$
satisfies the following condition: for any $c\in\varphi'(\mathcal{M},\bar{a})$,
$\models\chi(c,\bar{a},d)$ and $\chi(c,\bar{a},y)$ has finitely many, no more than $m=k_1+\ldots+k_n$
solutions. Consequently, the formula
$$\theta(\bar{a},y)\rightleftharpoons\exists
x(\chi(x,\bar{a},y)\wedge\forall
z((\varphi'(z,\bar{a})\to(\chi(x,\bar{a},y)\wedge\chi(z,\bar{a},y)))$$
satisfies $d$ and has no more than $m$ solutions. This fact contradicts the condition
$d\notin{\rm acl}(\bar{a}Z)$. \endproof

\medskip
Assuming that $b\notin{\rm acl}(aZ)$, we construct by induction a
countable model $\mathcal{N}\preccurlyeq\mathcal{M}$ such that
${\rm acl}(aZ)\subset N$, $b\notin N$, and
$N=\bigcup\limits_{n\in\omega}A_n$ for a chain of some sets $A_n$.

In the initial step we consider the set $A_0={\rm acl}(aZ)$ and
renumber all consistent formulas of the form $\varphi(x,\bar{a})$,
$\bar{a}\in A_0$: $\Phi_0=\{\varphi_{0,m}(x,\bar{a}_m)\mid
m\in\omega\}$. According this numeration construct at most a
countable set
$A_1=\bigcup\limits_{m\in\omega\cup\{-1\}}A_{1,m}\supset A_0$ with
the condition $b\notin{\rm acl}(A_1)$. Let
$A_{1,-1}\rightleftharpoons A_0$. If the set $A_{1,m-1}$ had been
already defined and $\varphi_{0,m}(\mathcal{M},\bar{a}_m)\cap
A_{1,m-1}\ne\varnothing$ then we put $A_{1,m}\rightleftharpoons
A_{1,m-1}$; if $\varphi_{0,m}(\mathcal{M},\bar{a}_m)\cap
A_{1,m-1}=\varnothing$ choose by Lemma \ref{le224} an element
$c_m\in\varphi_m(\mathcal{M},\bar{a}_m)$ such that $b\notin{\rm
acl}(c_mA_{1,m-1})$, and put $A_{1,m}\rightleftharpoons{\rm
acl}(c_mA_{1,m-1})$.

If at most a countable set $A_n$ had been already constructed,
renumber all consistent formulas of the form $\varphi(x,\bar{a})$,
$\bar{a}\in A_n$: $\Phi_n=\{\varphi_{n,m}(x,\bar{a}_m)\mid
m\in\omega\}$. According to this enumeration construct at most a
countable set
$A_{n+1}=\bigcup\limits_{m\in\omega\cup\{-1\}}A_{n+1,m}\supset
A_n$ with the condition $b\notin{\rm acl}(A_{n+1})$. Let
$A_{n+1,-1}\rightleftharpoons A_n$. If the set $A_{n+1,m-1}$ had
been already defined and $\varphi_{n,m}(\mathcal{M},\bar{a}_m)\cap
A_{n+1,m-1}\ne\varnothing$ then put $A_{n+1,m}\rightleftharpoons
A_{n+1,m-1}$; if $\varphi_{n,m}(\mathcal{M},\bar{a}_m)$ $\cap$
 $A_{n+1, m-1}=\varnothing$, choose by Lemma \ref{le224} an element
$c_m\in\varphi_{n,m}(\mathcal{M},\bar{a}_m)$ such that $b\notin{\rm acl}(c_mA_{n+1, m-1})$
and put $A_{n+1,m}\rightleftharpoons{\rm acl}(c_mA_{n+1,m-1})$.

By constructing the set $\bigcup\limits_{n\in\omega}A_n$ forms a required universe $N$ of a countable
model $\mathcal{N}\preccurlyeq\mathcal{M}$ such that ${\rm
acl}(Z\cup\{a\})\subseteq N$ and $b\notin N$. \endproof

\medskip
Applying Lemma \ref{le224}, we obtain the following lemma.

\begin{lemma}\label{le225}
Let $\mathcal{M}$ be an $\omega$-saturated model of a complete theory $T$, $\bar{a},\bar{b}\in M$,
$Z$ be the algebraic closure of some finite set in $\mathcal{M}$. If $({\rm
acl}(\bar{a}Z)\cap{\rm acl}(\bar{b}Z))\setminus Z=\varnothing$ and
$\varphi(x,\bar{a}')$ is a consistent formula, $\bar{a}'\in
\bar{a}Z$, then there is $c\in\varphi(\mathcal{M},\bar{a}')$
such that $({\rm acl}(\bar{a}cZ)\cap{\rm acl}(\bar{b}Z))\setminus
Z=\varnothing$.
\end{lemma}

\begin{theorem}\label{th226}
Let $\mathcal{M}$ be an $\omega$-saturated model of a countable complete theory $T$, $Z$ be the algebraic
closure of some finite set in $\mathcal{M}$, $a$ and $b$ be elements of $\mathcal{M}$, $a,b\notin Z$.
The following are equivalent:

$(1)$ the elements $a$ and $b$ are $Z$-separable in $H(\mathcal{M})$
by some sets $y$ and $z$ from $H(\mathcal{M})$ containing $Z$;

$(2)$ the elements $a$ and $b$ are $Z$-separable in $H_{\omega_1}(\mathcal{M})$ by some sets $y$ and
$z$ from $H_{\omega_1}(\mathcal{M})$ containing $Z$;

$(3)$ $({\rm acl}(aZ)\cap{\rm acl}(bZ))\setminus Z=\varnothing$.
\end{theorem}

Proof. As in the proof of Theorem \ref{th221} it is sufficiently
to prove the implication $(3)\Rightarrow(2)$. Assuming $({\rm
acl}(aZ)\cap{\rm acl}(bZ))\setminus Z=\varnothing$, we construct
by induction countable models
$\mathcal{N}_a,\mathcal{N}_b\preccurlyeq\mathcal{M}$ such that
${\rm acl}(aZ)\subseteq N_a$, ${\rm acl}(bZ)\subseteq N_b$,
$(N_a\cap N_b)\setminus Z=\varnothing$,
$N_a=\bigcup\limits_{n\in\omega}A_n$ for a chain of some sets
$A_n$ and $N_b=\bigcup\limits_{n\in\omega}B_n$ for a chain of some
sets $B_n$.

In the initial step we consider the sets $A_0={\rm acl}(aZ)$,
$B_0={\rm acl}(bZ)$ and enumerate all consistent formulas of the
form $\varphi(x,\bar{a})$, $\bar{a}\in A_0$:
$\Phi_0=\{\varphi_{0,m}(x,\bar{a}_m)\mid m\in\omega\}$. According
to this enumeration we construct at most countable set
$A_1=\bigcup\limits_{m\in\omega\cup\{-1\}}A_{1,m}\supset A_0$ with
the condition $({\rm acl}(A_1)\cap B_0)\setminus Z=\varnothing$.
Let $A_{1,-1}\rightleftharpoons A_0$. If the set $A_{1,m-1}$ had
been already defined and $\varphi_{0,m}(\mathcal{M},\bar{a}_m)\cap
A_{1,m-1}\ne\varnothing$, then put $A_{1,m}\rightleftharpoons
A_{1,m-1}$; if $\varphi_{0,m}(\mathcal{M},\bar{a}_m)\cap
A_{1,m-1}=\varnothing$ then by Lemma \ref{le225} we choose an
element $c_m\in\varphi_m(\mathcal{M},\bar{a}_m)$ such that $({\rm
acl}(c_mA_{1,m-1})\cap{\rm acl}(B_0))\setminus Z=\varnothing$ and
put $A_{1,m}\rightleftharpoons {\rm acl}(c_mA_{1,m-1})$.

If the set $A_1$ had been already defined, extend symmetrically the set $B_0$ to an algebraically closed
set $B_1$ such that $B_1\supseteq Z$, all consistent formulas $\varphi(x,\bar{b})$,
$\bar{b}\in B_0$, are realized in $B_1$ è $({\rm acl}(A_1)\cap{\rm
acl}(B_1))\setminus Z=\varnothing$.

If at most countable sets $A_n$ and $B_n$ had been already
constructed, renumber all consistent formulas of the form
$\varphi(x,\bar{a})$, $\bar{a}\in A_n$:
$\Phi_n=\{\varphi_{n,m}(x,\bar{a}_m)\mid m\in\omega\}$. According
to this numeration construct at most a countable set
$A_{n+1}=\bigcup\limits_{m\in\omega\cup\{-1\}}A_{n+1,m}\supset
A_n$ with the condition $({\rm acl}(A_{n+1})\cap{\rm
acl}(B_1))\setminus Z=\varnothing$. Let
$A_{n+1,-1}\rightleftharpoons A_n$. If the set $A_{n+1,m-1}$ had
been already defined and $\varphi_{0,m}(\mathcal{M},\bar{a}_m)\cap
A_{n+1,m-1}\ne\varnothing$, then put $A_{n+1,m}\rightleftharpoons
A_{n+1,m-1}$; if $\varphi_{0,m}(\mathcal{M},\bar{a}_m)\cap
A_{n+1,m-1}=\varnothing$, then by Lemma \ref{le225} choose an
element $c_m\in\varphi_{n,m}(\mathcal{M},\bar{a}_m)$ such that
$({\rm acl}(c_mA_{n+1,m-1})\cap{\rm acl}(B_n))\setminus
Z=\varnothing$, and put $A_{n+1,m}\rightleftharpoons
A_{n+1,m-1}\cup\{c_m\}$.

If we have the set $A_{n+1}$ then extend symmetrically the set $B_n$ to at most a countable set $B_{n+1}$
such that all consistent formulas $\varphi(x,\bar{b})$, $\bar{b}\in B_n$,
are realized in $B_{n+1}$ è  $({\rm acl}(A_{n+1})\cap{\rm
acl}(B_{n+1}))\setminus Z=\varnothing$.

By constructing the sets $\bigcup\limits_{n\in\omega}A_n$ and
$\bigcup\limits_{n\in\omega}B_n$ form required universes $N_a$ and $N_b$ respectively of $Z$-separable
countable models $\mathcal{N}_a,\mathcal{N}_b\preccurlyeq\mathcal{M}$ such that
$a\in N_a$ and $b\in N_b$. \endproof

Combining proofs of Claims \ref{co27}--\ref{co34} and Theorems \ref{th221}, \ref{th226}, we obtain the
following assertions.

\begin{corollary}\label{co47}
Let $\mathcal{M}$ be an $\omega$-saturated model of a small theory $T$,  $Z$ be the algebraic closure of some
finite set in $\mathcal{M}$, $a$ and $b$ be elements of $\mathcal{M}$, $a,b\notin Z$.
The following are equivalent:

$(1)$ the element $a$ is $Z$-separable from the element $b$ in $H(\mathcal{M})$
by some set $y$ from $H(\mathcal{M})$ containing $Z$;

$(2)$ the element $a$ is $Z$-separable from the element $b$ in $H_{\omega_1}(\mathcal{M})$
by some set $y$ from $H_{\omega_1}(\mathcal{M})$ containing $Z$;

$(3)$ the element $a$ is $Z$-separable from the element $b$ in $H_p(\mathcal{M})$
by some set $y$ from $H_p(\mathcal{M})$ containing $Z$;

$(4)$ $b\notin{\rm acl}(aZ)$.
\end{corollary}

\begin{corollary}\label{co48}
Let $\mathcal{M}$ be an $\omega$-saturated model of a small theory $T$, $Z$ be the algebraic closure of some finite
set in $\mathcal{M}$, $a$ and $b$ be elements of $\mathcal{M}$, $a,b\notin Z$. The following are equivalent:

$(1)$ the elements $a$ and $b$ are $Z$-separable in $H(\mathcal{M})$
by some sets $y$ and $z$ from $H(\mathcal{M})$ containing $Z$;

$(2)$ the elements $a$ and $b$ are $Z$-separable in $H_{\omega_1}(\mathcal{M})$ by some sets $y$ and
$z$ from $H_{\omega_1}(\mathcal{M})$ containing $Z$;

$(3)$ the elements $a$ and $b$ are separable in $H_p(\mathcal{M})$ by some sets $y$ and $z$ from
$H_p(\mathcal{M})$ containing $Z$;

$(4)$ $({\rm acl}(aZ)\cap{\rm acl}(bZ))\setminus Z=\varnothing$.
\end{corollary}

\begin{theorem}\label{th51} Let $\mathcal{M}$ be a $\lambda$-saturated model of a complete theory $T$,
$\lambda\geq{\rm max}\{|\Sigma(T)|$, $\omega\}$, $A$ and $B$ be
nonempty sets in $\mathcal{M}$ having cardinalities $<\lambda$,
$Z$ be the algebraic closure of some set of cardinality $<\lambda$
in $\mathcal{M}$. The following are equivalent:

$(1)$ the set $A$ is $Z$-separable from the set $B$ in $H(\mathcal{M})$ by some set $y$ from
$H(\mathcal{M})$ containing $Z$;

$(2)$ the set $A$ is $Z$-separable from the set $B$ in $H_\lambda(\mathcal{M})$ by some set $y$ from
$H_\lambda(\mathcal{M})$ containing $Z$;

$(3)$ $({\rm acl}(A\cup Z)\cap B)\setminus Z=\varnothing$.
\end{theorem}

\begin{theorem}\label{th52} Let $\mathcal{M}$ be a $\lambda$-saturated model of a complete theory $T$,
$\lambda\geq{\rm max}\{|\Sigma(T)|$, $\omega\}$, $A$ and $B$ be
nonempty sets in $\mathcal{M}$ having cardinalities $<\lambda$,
$Z$ be the algebraic closure of some set of cardinality $<\lambda$
in $\mathcal{M}$. The following are equivalent:

$(1)$ the sets $A$ and $B$ are $Z$-separable in $H(\mathcal{M})$ by some sets $y$ and $z$ from
$H(\mathcal{M})$ containing $Z$;

$(2)$ the sets $A$ and $B$ are $Z$-separable in $H_\lambda(\mathcal{M})$ by some sets $y$ and
$z$ from $H_\lambda(\mathcal{M})$ containing $Z$;

$(3)$ $({\rm acl}(A\cup Z)\cap{\rm acl}(B\cup Z))\setminus
Z=\varnothing$.
\end{theorem}

\begin{corollary}\label{co53}
Let $\mathcal{M}$ be an $\omega$-saturated model of a small theory
$T$, $A$ and $B$ be finite nonempty sets in $\mathcal{M}$, $Z$ be
the algebraic closure of some finite set in $\mathcal{M}$. The
following are equivalent:

$(1)$ the set $A$ is $Z$-separable from the set $B$ in $H(\mathcal{M})$ by some set $y$ from
$H(\mathcal{M})$ containing $Z$;

$(2)$ the set $A$ is $Z$-separable from the set $B$ in $H_{\omega_1}(\mathcal{M})$ by some set $y$
from $H_{\omega_1}(\mathcal{M})$ containing $Z$;

$(3)$ the set $A$ is $Z$-separable from the set $B$ in $H_p(\mathcal{M})$ by some set $y$ from
$H_p(\mathcal{M})$ containing $Z$;

$(4)$ $({\rm acl}(A\cup Z)\cap B)\setminus Z=\varnothing$.
\end{corollary}

\begin{corollary}\label{co54}
Let $\mathcal{M}$ be an $\omega$-saturated model of a small theory
$T$, $A$ and $B$ be finite nonempty sets in $\mathcal{M}$, $Z$ be
the algebraic closure of some finite set in $\mathcal{M}$. The
following are equivalent:

$(1)$ the sets $A$ and $B$ are $Z$-separable in $H(\mathcal{M})$ by some sets $y$ and $z$ from
$H(\mathcal{M})$ containing $Z$;

$(2)$ the sets $A$ and $B$ are $Z$-separable in $H_{\omega_1}(\mathcal{M})$ by some sets $y$ and
$z$ from $H_{\omega_1}(\mathcal{M})$ containing $Z$;

$(3)$ the sets $A$ and $B$ are $Z$-separable in $H_p(\mathcal{M})$ by some sets $y$ and $z$
from $H_p(\mathcal{M})$ containing $Z$;

$(4)$ $({\rm acl}(A\cup Z)\cap{\rm acl}(B\cup Z))\setminus
Z=\varnothing$.
\end{corollary}

\section{On separability in hypergraphs of models of ordered theories}

\begin{definition} \rm  \cite{CCMCT14, IPT}  Let $p_1(x_1),\ldots,p_n(x_n)\in S_1(T)$.
A type $q(x_1,\ldots,x_n)\in S(T)$ is called \emph{$(p_1,\ldots,p_n)$-type} if
$q(x_1,\ldots,x_n)\supseteq\bigcup\limits_{i=1}^n p_i(x_i)$. The set of all
\ $(p_1,\ldots,p_n)$-types of a theory  $T$ is denoted by $S_{p_1,\ldots,p_n}(T)$.
A countable theory $T$ is called {\em almost $\omega$-categorical} if for any types
$p_1(x_1),\ldots,p_n(x_n)\in S(T)$ there exist only finitely many types
$q(x_1,\ldots,x_n)\in S_{p_1,\ldots,p_n}(T)$.
\end{definition}

\begin{theorem} \label{th_sep} Let $T$ be an almost $\omega$-categorical quite o-minimal
theory, $\mathcal{M}$ be an $\omega$-saturated model of the theory $T$, $Z$ be the algebraic
closure of some finite set in $\mathcal{M}$, $a,b\in M\setminus Z$.
Then the following conditions are equivalent:

$(1)$ $a$ is $Z$-separable from $b$ in $\mathcal{H}(\mathcal{M})$
by some set $y$ from $H(\mathcal{M})$ containing $Z$;

$(2)$ $b$ is $Z$-separable from $a$ in $\mathcal{H}(\mathcal{M})$
by some set $y$ from $H(\mathcal{M})$ containing $Z$;

$(3)$ the elements $a$ and $b$ are $Z$-separable in $H(\mathcal{M})$
by some sets $y$ and $z$ from $H(\mathcal{M})$ containing $Z$;

$(4)$ $a\not\in {\rm dcl}(\{bZ\})$;

$(5)$ $b\not\in {\rm dcl}(\{aZ\})$.

$(6)$ $({\rm dcl}(aZ)\cap{\rm dcl}(bZ))\setminus Z=\varnothing$.

\end{theorem}

Proof of Theorem \ref{th_sep}. By Proposition 3.9
\cite{ks_mn_2017} Exchange Principle for algebraic closure holds. By linear ordering of
the model $\mathcal{M}$ ${\rm dcl}(A)={\rm acl}(A)$ for any $A\subseteq M$. Then by proofs
of Theorems \ref{th221} and \ref{th226} we have an equivalence of the conditions (1)--(6).
\endproof

\begin{note}\label{not71} {\rm 1. Theorem \ref{th_sep} remains true for an arbitrary theory
satisfying both Exchange Principle for algebraic closure and the condition
${\rm dcl}(A)={\rm acl}(A)$ for any $A\subseteq M$.

2. If Exchange Principle for algebraic closure holds and the condition ${\rm
dcl}(A)={\rm acl}(A)$ for any $A\subseteq M$ doesn't hold, Theorem
\ref{th_sep} remains true if we replace ${\rm dcl}$ by ${\rm acl}$.

3. If the condition ${\rm dcl}(A)={\rm acl}(A)$ for any $A\subseteq M$ holds and Exchange
Principle for algebraic closure doesn't hold, Theorem \ref{th_sep} splits into three
independent statements $(1)\Leftrightarrow(5)$, $(2)\Leftrightarrow(4)$,
$(3)\Leftrightarrow(6)$.}
\end{note}

Theorem \ref{th_sep} immediately implies the following

\begin{corollary} \label{cor_sep} Let $T$ be an almost $\omega$-categorical quite o-minimal
theory, $\mathcal{M}$ be an $\omega$-saturated model of the theory $T$, $a,b\in M\setminus dcl(\varnothing)$.
Then the following conditions are equivalent:

$(1)$ $a$ is separable from $b$ in $\mathcal{H}(\mathcal{M})$;

$(2)$ $b$ is separable from $a$ in $\mathcal{H}(\mathcal{M})$;

$(3)$ $a\not\in \dcl(\{b\})$;

$(4)$ $b\not\in \dcl(\{a\})$.

\end{corollary}

\begin{example}\label{ex1}\rm \cite{mms}
Let $\mathcal{M}=\langle M;<, P^1_1, P^1_2, f^1\rangle$ be a linearly ordered structure such that
$M$ is the disjoint union of interpretations of unary predicates $P_1$ and $P_2$ so that
$P_1(\mathcal{M})<P_2(\mathcal{M})$. We identify an interpretation of $P_2$ with the set of rational
numbers $\mathbb{Q}$, ordered as usual, and $P_1$ with $\mathbb{Q}\times \mathbb{Q}$, ordered
lexicographically. The symbol $f$ is interpreted by a partial unary function with
$Dom(f)=P_1(\mathcal{M})$ and $Range(f)=P_2(\mathcal{M})$ and is defined by the equality
$f((n,m))=n$ for all $(n,m)\in \mathbb{Q}\times \mathbb{Q}$.
\end{example}

It is known that $\mathcal{M}$ is a countably categorical weakly
o-minimal structure, and $Th(\mathcal{M})$ is not quite o-minimal.
Take arbitrary $a \in P_1(\mathcal{M})$, $b\in P_2(\mathcal{M})$
such that $f(a)=b$. Then we obtain that $a$ is separable from $b$
in $\mathcal{H(\mathcal{M})}$, but $b$ is not separable from $a$
in $\mathcal{H(\mathcal{M})}$.

\begin{proposition}\label{pr_fsets}
Let $T$ be an almost $\omega$-categorical quite o-minimal theory,
$\mathcal{M}\models T$, $A=\{a_1, \ldots, a_{n_1}\}$, $B=\{b_1,
\ldots, b_{n_2}\}\subseteq M$ for some positive $n_1, n_2<\omega$.
Then the following conditions are equivalent:

$(1)$ $A$ and $B$ are separable from each other in $\mathcal{H}(\mathcal{M})$;

$(2)$ $\dcl(A)\cap \dcl(B)=\emptyset$.

$(3)$ $\dcl(\{a_i\})\cap\dcl(\{b_j\})=\emptyset$ for any $1\le i\le n_1$, $1\le j\le n_2$.
\end{proposition}

Proof of Proposition \ref{pr_fsets}. $(1)\Rightarrow (2)$ Let $A$
be separable from $B$ in $\mathcal{H}(\mathcal{M})$. This means
that there is $\mathcal{M}_1\prec \mathcal{M}$ such that
$A\subseteq M_1$ and $B\cap M_1=\emptyset$. Then we have:
$\dcl(A)\subseteq M_1$, whence we obtain that $\dcl(A)\cap
B=\emptyset$. Similarly, by the condition of separability of $B$
from $A$ in $\mathcal{H}(\mathcal{M})$ it can be established that
$\dcl(B)\cap A=\emptyset$.

Assume the contrary: $\dcl(A)\cap \dcl(B)\ne \emptyset$.
Consequently, there is $c\in M$ such that $c\in \dcl(A)$ and
$c\in \dcl(B)$. But then by binarity of $Th(\mathcal{M})$
there exist $a\in A$ and $b\in B$ such that $c\in \dcl(\{a\})$ and $c\in \dcl(\{b\})$.
By holding Exchange Principle for algebraic closure we obtain that $b\in \dcl(\{a\})$.
The last contradicts the condition $\dcl(A)\cap B=\emptyset$.

$(2)\Rightarrow (1)$ In this case we assert that $M_1:=M\setminus
\dcl(A)$ and $M_2:=M\setminus \dcl(B)$ are universes of elementary submodels of the model $\mathcal{M}$.

$(2)\Leftrightarrow (3)$ By binarity of $Th(\mathcal{M})$. \endproof

\begin{proposition}\label{pr_hfsets}
Let $T$ be an almost $\omega$-categorical quite o-minimal theory,
$\mathcal{M}\models T$, $Z=\dcl(\emptyset)$, $A=\{a_1, \ldots,
a_{n_1}\}$, $B=\{b_1, \ldots, b_{n_2}\}\subseteq M$ for some
positive $n_1, n_2<\omega$ so that $A\cap Z=B\cap Z=\emptyset$.
Then the following conditions are equivalent:

$(1)$ $A$ and $B$ are $Z$-separable in $\mathcal{H}(\mathcal{M})$;

$(2)$ $\dcl(A)\cap \dcl(B)=Z$.

$(3)$ $\dcl(\{a_i\})\cap\dcl(\{b_j\})=Z$ for any $1\le i\le n_1$, $1\le j\le n_2$.
\end{proposition}

Proof of Proposition \ref{pr_hfsets}. $(1)\Rightarrow (2)$ Let $A$ and $B$ be $Z$-separable in
$\mathcal{H}(\mathcal{M})$. Then there exist $\mathcal{M}_1, \mathcal{M}_2\prec \mathcal{M}$ such
that $(\mathcal{M}_1\cap\mathcal{M}_2)\setminus Z=\emptyset$, $A\subseteq M_1$ and $B\subseteq M_2$.
Consequently, $\dcl(A)\cap \dcl(B)\subseteq \mathcal{M}_1\cap \mathcal{M}_2$. Then
$[\dcl(A)\cap \dcl(B)]\setminus Z=\emptyset$, whence $\dcl(A)\cap \dcl(B)=Z$.

$(2)\Rightarrow (1)$ In this case we assert that $M_1:=[M\setminus \dcl(A)]\cup Z$ and
$M_2:=[M\setminus \dcl(B)]\cup Z$ are universes of elementary submodels of the model $\mathcal{M}$. \endproof

\medskip
Arguments for Propositions \ref{genpr_asets} and \ref{pr_fsets}
imply the following

\begin{proposition}\label{pr_asets}
Let $T$ be an almost $\omega$-categorical quite o-minimal theory,
$\mathcal{M}\models T$, $\emptyset\ne A,B\subseteq M$, $\mathcal{M}$ be $|A\cup
B|^+$-saturated. Then $A$ and $B$ are separable from each other in
$\mathcal{H}(\mathcal{M})$ if and only if the following conditions
hold:

$(1)$ $\dcl(\{a\})\cap \dcl(\{b\})=\emptyset$ for any $a\in A$ and $b\in B$;

$(2)$ For any isolated type $p\in S_1(\emptyset)$,
$p(\mathcal{M})\setminus \dcl(A)\ne \emptyset$ and
$p(\mathcal{M})\setminus \dcl(B)\ne \emptyset$.
\end{proposition}

\begin{corollary}\label{pr_ahsets}
Let $T$ be an almost $\omega$-categorical quite o-minimal theory,
$\mathcal{M}\models T$, $Z=\dcl(\emptyset)$, $A,B$ be non-empty subsets of $M$
such that $A\cap Z=B\cap Z=\emptyset$, $\mathcal{M}$ be $|A\cup
B|^+$-saturated. Then $A$ and $B$ are $Z$-separable in
$\mathcal{H}(\mathcal{M})$ if and only if the following conditions
hold:

$(1)$ $\dcl(\{a\})\cap \dcl(\{b\})=Z$ for any $a\in A$ and $b\in B$;

$(2)$ For any isolated type $p\in S_1(\emptyset)$,
$p(\mathcal{M})\setminus \dcl(A)\ne \emptyset$ and
$p(\mathcal{M})\setminus \dcl(B)\ne \emptyset$.
\end{corollary}

Arguments for Propositions \ref{genpr_asets} and \ref{pr_fsets} as
well as Theorem \ref{th52} imply the following

\begin{proposition}\label{pr_Z_sep}
Let $T$ be an almost $\omega$-categorical quite o-minimal theory,
$\mathcal{M}\models T$ be $\lambda$-saturated, $\lambda\geq{\rm
max}\{|\Sigma(T)|,\omega\}$, $A$ and $B$ be nonempty sets in
$\mathcal{M}$ having cardinalities $<\lambda$, $Z$ be the
algebraic closure of some set of cardinality $<\lambda$ in
$\mathcal{M}$. Then the following are equivalent:

$(1)$ $A$ and $B$ are $Z$-separable in $\mathcal{H}(\mathcal{M})$;

$(2)$ $(\dcl(aZ)\cap \dcl(bZ))\setminus Z=\emptyset$ for any $a\in
A$ and $b\in B$.
\end{proposition}

\vskip 2mm {\bf Acknowledgements.} This research was partially
supported by Committee of Science in Education and Science
Ministry of the Republic of Kazakhstan (Grant No. AP05132546) and
Russian Foundation for Basic Researches (Project No.
17-01-00531-a).

\end{document}